\documentstyle[12pt]{article}
\begin{document}
\newcommand{\text}[1]{\mbox{{\rm #1}}}
\newcommand{\gd}{\delta}
\newcommand{\itms}[1]{\item[[#1]]}
\newcommand{\nin}{\in\!\!\!\!\!/}
\newcommand{\sub}{\subset}
\newcommand{\cntd}{\subseteq}
\newcommand{\go}{\omega}
\newcommand{\Pa}{P_{a^\nu,1}(U)}
\newcommand{\fx}{f(x)}
\newcommand{\fy}{f(y)}
\newcommand{\gD}{\Delta}
\newcommand{\gl}{\lambda}
\newcommand{\gL}{\Lambda}
\newcommand{\half}{\frac{1}{2}}
\newcommand{\sto}[1]{#1^{(1)}}
\newcommand{\stt}[1]{#1^{(2)}}
\newcommand{\Z}{\hbox{\sf Z\kern-0.720em\hbox{ Z}}}
\newcommand{\singcolb}[2]{\left(\begin{array}{c}#1\\#2
\end{array}\right)}
\newcommand{\ga}{\alpha}
\newcommand{\gb}{\beta}
\newcommand{\gga}{\gamma}
\newcommand{\ul}{\underline}
\newcommand{\ol}{\overline}
\newcommand{\qed}{\kern 5pt\vrule height8pt width6.5pt depth2pt}
\newcommand{\Lrraro}{\Longrightarrow}
\newcommand{\Nb}{|\!\!/}
\newcommand{\NN}{{\rm I\!N}}
\newcommand{\bsl}{\backslash}
\newcommand{\gt}{\theta}
\newcommand{\op}{\oplus}
\newcommand{\C}{{\bf C}}
\newcommand{\Q}{{\bf Q}}
\newcommand{\Op}{\bigoplus}
\newcommand{\CR}{{\cal R}}
\newcommand{\tr}{\bigtriangleup}
\newcommand{\grr}{\omega_1}
\newcommand{\ben}{\begin{enumerate}}
\newcommand{\een}{\end{enumerate}}
\newcommand{\ndiv}{\not\mid}
\newcommand{\bab}{\bowtie}
\newcommand{\hal}{\leftharpoonup}
\newcommand{\har}{\rightharpoonup}
\newcommand{\ot}{\otimes}
\newcommand{\OT}{\bigotimes}
\newcommand{\bwe}{\bigwedge}
\newcommand{\gep}{\varepsilon}
\newcommand{\gs}{\sigma}
\newcommand{\rbraces}[1]{\left( #1 \right)}
\newcommand{\bbox}{$\;\;\rule{2mm}{2mm}$}
\newcommand{\sbraces}[1]{\left[ #1 \right]}
\newcommand{\bbraces}[1]{\left\{ #1 \right\}}
\newcommand{\OO}{_{(1)}}
\newcommand{\TT}{_{(2)}}
\newcommand{\FF}{_{(3)}}
\newcommand{\minus}{^{-1}}
\newcommand{\CV}{\cal V}
\newcommand{\CVs}{\cal{V}_s}
\newcommand{\un}{U_q(sl_n)'}
\newcommand{\on}{O_q(SL_n)'}
\newcommand{\slq}{U_q(sl_2)}
\newcommand{\olq}{O_q(SL_2)}
\newcommand{\UU}{U_{(N,\nu,\go)}}
\newcommand{\HH}{H_{n,q,N,\nu}}
\newcommand{\ct}{\centerline}
\newcommand{\bs}{\bigskip}
\newcommand{\qua}{\rm quasitriangular}
\newcommand{\ms}{\medskip}
\newcommand{\noin}{\noindent}
\newcommand{\mat}[1]{$\;{#1}\;$}
\newcommand{\raro}{\rightarrow}
\newcommand{\map}[3]{{#1}\::\:{#2}\raro{#3}}
\newcommand{\alg}{{\rm Alg}}
\def\newtheorems{\newtheorem{theorem}{Theorem}[subsection]
                 \newtheorem{cor}[theorem]{Corollary}
                 \newtheorem{prop}[theorem]{Proposition}
                 \newtheorem{lemma}[theorem]{Lemma}
                 \newtheorem{defn}[theorem]{Definition}
                 \newtheorem{Theorem}{Theorem}[section]
                 \newtheorem{Corollary}[Theorem]{Corollary}
                 \newtheorem{Proposition}[Theorem]{Proposition}
                 \newtheorem{Lemma}[Theorem]{Lemma}
                 \newtheorem{Defn}[Theorem]{Definition}
                 \newtheorem{Example}[Theorem]{Example}
                 \newtheorem{Remark}[Theorem]{Remark}
                 \newtheorem{claim}[theorem]{Claim}
                 \newtheorem{sublemma}[theorem]{Sublemma}
                 \newtheorem{example}[theorem]{Example}
                 \newtheorem{remark}[theorem]{Remark}
                 \newtheorem{question}[theorem]{Question}
                 \newtheorem{Question}[Theorem]{Question}
                 \newtheorem{conjecture}{Conjecture}[subsection]}
\newtheorems
\newcommand{\proof}{\par\noindent{\bf Proof:}\quad}
\newcommand{\dmatr}[2]{\left(\begin{array}{c}{#1}\\
                            {#2}\end{array}\right)}
\newcommand{\doubcolb}[4]{\left(\begin{array}{cc}#1&#2\\
#3&#4\end{array}\right)}
\newcommand{\qmatrl}[4]{\left(\begin{array}{ll}{#1}&{#2}\\
                            {#3}&{#4}\end{array}\right)}
\newcommand{\qmatrc}[4]{\left(\begin{array}{cc}{#1}&{#2}\\
                            {#3}&{#4}\end{array}\right)}
\newcommand{\qmatrr}[4]{\left(\begin{array}{rr}{#1}&{#2}\\
                            {#3}&{#4}\end{array}\right)}
\newcommand{\smatr}[2]{\left(\begin{array}{c}{#1}\\
                            \vdots\\{#2}\end{array}\right)}

\newcommand{\ddet}[2]{\left[\begin{array}{c}{#1}\\
                           {#2}\end{array}\right]}
\newcommand{\qdetl}[4]{\left[\begin{array}{ll}{#1}&{#2}\\
                           {#3}&{#4}\end{array}\right]}
\newcommand{\qdetc}[4]{\left[\begin{array}{cc}{#1}&{#2}\\
                           {#3}&{#4}\end{array}\right]}
\newcommand{\qdetr}[4]{\left[\begin{array}{rr}{#1}&{#2}\\
                           {#3}&{#4}\end{array}\right]}

\newcommand{\qbracl}[4]{\left\{\begin{array}{ll}{#1}&{#2}\\
                           {#3}&{#4}\end{array}\right.}
\newcommand{\qbracr}[4]{\left.\begin{array}{ll}{#1}&{#2}\\
                           {#3}&{#4}\end{array}\right\}}

\title{The Classification of Triangular Semisimple and Cosemisimple Hopf
Algebras Over an Algebraically Closed Field}
\author{Pavel Etingof
\\Department of Mathematics, Rm 2-165\\
Massachusetts Institute of Technology\\
77 Massachusetts Avenue\\
Cambridge, MA 02139
\and Shlomo Gelaki $^1$\\Mathematical Sciences Research Institute\\
1000 Centennial Drive\\Berkeley, CA 94720}
\footnotetext[1]{This research was supported by THE ISRAEL SCIENCE
FOUNDATION founded by the Israel Academy of Sciences and Humanities.} 
\date{November 2, 1999}
\maketitle

\section{Introduction}
In this paper we classify triangular semisimple
and cosemisimple Hopf algebras over
{\em any} algebraically closed field $k$. Namely, we construct,
for each positive integer $N$, relatively
prime to the characteristic of $k$ if it is positive, a bijection between
the set of
isomorphism classes of triangular semisimple and cosemisimple Hopf algebras
of dimension $N$ over $k$, and the set of isomorphism classes
of quadruples $(G,H,V,u)$, where $G$ is a group of order $N$, $H$ is a
subgroup of $G$, $V$ is an irreducible projective representation of $H$
over $k$ of dimension $|H|^{1/2}$, and $u\in G$ is a central element of
order $\le 2.$ This classification
implies, in particular, that {\em any} triangular semisimple and
cosemisimple Hopf
algebra over $k$ can be obtained from a group algebra by a twist (it was
previously known only in characteristic $0$ [EG1, Theorem 2.1]).
It also implies that $(G,H,V,u)$ corresponds to a {\em
minimal} triangular 
semisimple Hopf algebra over $k$ if and only if $G$ is generated 
by $H$ and $u.$ We then answer positively the question from
[EG2] whether 
the group underlying a {\em minimal} triangular semisimple Hopf
algebra is solvable by proving that the group $H$ is a quotient 
of a central type group and hence solvable. 
We then conclude that any triangular semisimple
and cosemisimple Hopf algebra over $k$ of dimension bigger than $1$
contains a non-trivial grouplike element.

The classification uses Deligne's theorem on Tannakian
categories [De]
and the results of the paper [M] in an essential way. The proof of
solvability and existence of grouplike elements relies on a
theorem from [HI], which is proved using the classification
of finite simple groups. The classification in positive characteristic
relies also on the lifting functor from [EG4].

Throughout the paper, the ground field $k$ is
assumed to be algebraically closed.

{\bf Acknowledgments.} We are grateful to Donald Passman for
telling us about central type groups and the paper [HI]. We thank
Mikhail Movshev for useful comments on his results, and
Robert Guralnick
for his interest and useful conversations.
\section{Twists}
Let $A$ be a Hopf algebra over a field $k.$ Recall [Dr1]
that a {\em twist} for $A$ is an
invertible element $J\in A\ot A$ which satisfies
\begin{equation}\label{1}
(\Delta\ot I)(J)J_{12}=(I\ot \Delta)(J)J_{23}\;\; and \;\; (\varepsilon\ot
I)(J)=(I\ot \varepsilon)(J)=1,
\end{equation}
where $I$ is the identity map of $A.$

If $J$ is a twist for $A$ and $x$ is an invertible element of $A$ then
$J^x=\Delta(x)J(x^{-1}\otimes x^{-1})$ is also a twist for $A.$
We will call the twists $J$ and $J^x$ {\em gauge equivalent}.
The element $x$ will be called a {\em gauge transformation}.

Given a twist $J$ for $A$, one can define a Hopf algebra
$(A^J,\Delta^J,\varepsilon)$ as follows: $A^J=A$ as an
algebra, the coproduct is determined by
$$\Delta^J(x)=J^{-1}\Delta(x)J$$ for all $x\in A,$ and $\varepsilon$ is
the
ordinary counit of $A$. If $A$ is triangular with the
universal $R-$matrix $R$, then so is $A^J$, with
the universal $R-$matrix $R^J=J_{21}^{-1}RJ$. It is obvious that two
gauge
equivalent twists, when applied to a fixed
(triangular) Hopf algebra, produce two isomorphic (triangular) Hopf algebras.

Let $A=k[H]$ be the group algebra of a finite group $H.$
We will say that a twist
$J$ for $A$ is {\em minimal} if the right (and left) components of
the $R-$matrix $R^J=J_{21}^{-1}J$
span $A,$ i.e. if the corresponding triangular Hopf algebra
$(A^J,J_{21}^{-1}J)$ is minimal [R].

Let $(A,R)$ be any triangular semisimple and cosemisimple Hopf algebra
over $k$. Then the Drinfeld element $u$ of $A$ is a
grouplike element of order $\le 2$. Moreover, by [LR] in characteristic
$0,$ and by [EG4, Theorem 3.1] in positive characteristic, the square of
the antipode of $A$ is the identity map, and hence $u$ is
central.

If $k$ does not have characteristic $2,$ set 
$$R_u=\frac{1}{2}(1\otimes
1+1\otimes u+u\otimes 1-u\otimes u).$$ If $k$ is of
characteristic $2$ (in which case $u=1$ by semisimplicity),
set $R_u=1$. Then $(A,RR_u)$ is triangular semisimple and
cosemisimple with Drinfeld element $1.$ This observation
allows to reduce questions about triangular semisimple
and cosemisimple Hopf algebras over $k$ to the case when the
Drinfeld element is $1.$
\section{Twists for Group Algebras}
In this section we will prove the following theorem, which will be used
later for classification.
\begin{Theorem}\label{uni}
Let $G,$ $G'$ be finite groups,
$H,$ $H'$ subgroups of $G,$ $G'$ respectively, and
$J,$ $J'$ minimal twists for $k[H],$ $k[H']$ respectively.
Suppose that the triangular Hopf algebras
$k[G]^J,$ $k[G']^{J'}$ are isomorphic.
Then there exists a group isomorphism $\phi:G\to G'$ such that
$\phi(H)=H'$, and $(\phi\ot \phi)(J)$ is gauge equivalent to $J'$
as twists for $k[H']$.
\end{Theorem}

The rest of the section is devoted to the proof of the theorem.
\begin{Lemma}\label{unifunctor}
Let $\cal C$ be the category of $k-$representations
of a finite group, and $F_1,$ $F_2:{\cal C}\raro Vect(k)$ be two
fiber functors (i.e. exact and 
faithful symmetric tensor functors,
see [DM]) from $\cal C$ to the category of $k-$vector spaces. Then $F_1$
is isomorphic to $F_2.$
\end{Lemma}
\proof This is a special case of [DM, Theorem 3.2]. \qed

The following corollary of this lemma answers positively Movshev's
question [M, Remark 1] whether any symmetric twist is trivial.
\begin{Corollary}\label{symtwist} Let $G$ be a finite group, and $J$ be a
symmetric twist for $k[G]$ (i.e. $J_{21}=J$). Then $J$ is gauge equivalent
to $1\ot 1.$
\end{Corollary}
\proof Let $\cal C$ be the category of representations of $G.$
We have two symmetric tensor structures on the forgetful functor
$F:{\cal C}\raro Vect(k)$; namely, the trivial one and the one defined by
$J.$ By Lemma \ref{unifunctor}, the two fiber functors corresponding
to these structures are isomorphic.
But by definition, an isomorphism between them is
an invertible element $x\in k[G]$ such that $J=\Delta(x)(x^{-1}\ot
x^{-1}).$ \qed
\begin{Remark} {\rm Here is another proof of Corollary \ref{symtwist}
(in the case when the characteristic of $k$ is relatively prime to $|G|$)
which does
not use Lemma \ref{unifunctor} but uses the results of [M].
Consider the $G$-coalgebra $B_J=k[G]$ with coproduct $\tilde\Delta(x)=
(x\ot x)J$, and the dual algebra $B_J^*$.
According to [M], this algebra is semisimple,
$G$ acts transitively on its simple ideals, and $B_J^*$, along with the action
of $G,$ completely determines $J$ up to gauge transformations. Clearly,
since $J$ is symmetric, this algebra is commutative. So, it is
isomorphic, as a $G$-algebra, to the algebra of functions
on a set $X$ on which $G$ acts simply transitively. Corollary
\ref{symtwist} now
follows from the fact that such a $G$-set is unique up to an isomorphism
(the group $G$ itself with $G$ acting by left multiplication).}\qed
\end{Remark}
\begin{Lemma}\label{uniqueness}
Let $G,$ $G'$ be finite groups, $J,$ $J'$ twists for $k[G],$ $k[G']$
respectively, and suppose that the triangular Hopf algebras
$k[G]^J,$ $k[G']^{J'}$ are isomorphic.
Then there exists a group isomorphism $\phi:G\to G'$ such that
$(\phi\ot \phi)(J)$ is gauge equivalent to $J'.$
\end{Lemma}
\proof Let $f:k[G]^J\to k[G']^{J'}$ be an isomorphism
of triangular Hopf algebras.
Then $f$ defines an isomorphism of triangular Hopf algebras
$k[G]\to k[G']^{J'(f\ot f)(J)^{-1}}$. This implies that
the element $J'(f\ot f)(J)^{-1}$ is a symmetric twist
for $k[G']$. Thus, for some invertible $x\in k[G']$ one has
$J'(f\ot f)(J)^{-1}=\Delta(x)(x^{-1}\ot x^{-1})$.
Let $\phi=Ad(x^{-1})\circ f:k[G]\to k[G']$.
It is obvious that $\phi$ is a Hopf algebra isomorphism, so it comes
from a group isomorphism $\phi:G\to G'$.
We have $(\phi\ot\phi)(J)=\Delta(x)^{-1}J'(x\ot x)$, as desired. \qed

We can now prove Theorem \ref{uni}. By Lemma \ref{uniqueness},
it is sufficient to assume that $G'=G$, and that $J$ is gauge equivalent
to $J'$ as twists for $k[G],$ and it is enough to show that there exists
an element $a\in G$ such that $aHa^{-1}=H'$ and $(a\ot a)J(a^{-1}\ot a^{-1})$
is gauge equivalent to $J'$ as twists for $k[H'].$

So let $x\in k[G]$ be the invertible element such that
$\Delta(x)J(x^{-1}\ot x^{-1})=J'.$
In particular, this implies that $(x\ot x)R(x^{-1}\ot x^{-1})=R',$
where $R,R'$ are the $R$-matrices corresponding to $J,J'$ respectively. By
the
minimality of $J,$ $J',$ we have
$xk[H]x^{-1}=k[H']$. Thus, $$J_0=\Delta(x)(x^{-1}\ot x^{-1})=
J'(x\ot x)J^{-1}(x^{-1}\ot x^{-1})\in k[H']^{\ot 2}.$$
It is obvious that $J_0$ is a symmetric twist for $k[H']$, so by Corollary
\ref{symtwist},
it is gauge equivalent to $1\ot 1.$ Thus, $x=x_0a$, for some 
invertible $x_0\in k[H']$,
and $a\in G$. It is clear that $aHa^{-1}=H'$, and
$\Delta(x_0^{-1})J'(x_0\ot x_0)=(a\ot a)J(a^{-1}\ot a^{-1}).$
This concludes the proof of Theorem \ref{uni}.
\section{Construction of Triangular Semisimple and Cosemisi-mple Hopf
Algebras from Group-Theoretical Data}
Let $H$ be a finite group
such that $|H|$ is not divisible by the characteristic of $k.$
Suppose that $V$ is an irreducible projective representation of $H$ over
$k$ satisfying $dimV=|H|^{1/2}.$ Following the idea of the proof of
Proposition 5 in [M], we will construct a twist $J\in k[H]\ot k[H].$

Let $\pi :H\raro PGL(V)$ be the projective action of $H$ on $V$, and let
$\tilde \pi:H\to SL(V)$ be any lifting of this action
($\tilde\pi$ need not be a homomorphism).
We have $\tilde\pi(x)\tilde\pi(y)=c(x,y)\tilde\pi(xy)$, where $c$ is a
2-cocycle with coefficients in $k^*.$ By [M, Proposition 11], this
cocycle is nondegenerate
(see [M]) and hence by [M, Proposition 12] the representation
of $H$ on $End_{k}(V)$ is
isomorphic to the regular representation.
\begin{Remark} {\rm In the paper [M] it is assumed that the
characteristic of $k$ is equal to $0,$ but all the results generalize in
a straightforward way to the case when the characteristic of $k$ is
positive and relatively prime to the order of the group.}
\end{Remark}

Consider the simple coalgebra $(End_{k}V)^*$ with comultiplication
$\Delta.$ Clearly $H$ acts on
this coalgebra, and this representation of $H$ is also isomorphic to the
regular representation.
In particular we can choose an element $\lambda\in (End_{k}V)^*$ such
that the set $\{a\cdot \lambda|a\in H\}$ forms a basis of
$(End_{k}V)^*,$ and
$<\lambda,I>=1$ where $I$ is the unit element of $End_{k}V.$ Now, write
$\Delta(\lambda)=\sum_{a,b\in H} \gamma(a,b) a\cdot \lambda\ot b\cdot
\lambda,$ and set $J=\sum_{a,b\in H}
\gamma(a,b) a\ot b\in k[H]\ot k[H].$ We claim that $J$ is a
twist for $k[H]$.

Indeed, let $\tilde{\Delta}:k[H]\raro k[H]^{\ot 2}$ be determined by 
$a\mapsto
(a\ot a)J,$ and let $f:k[H]\raro (End_{k}V)^*$ be determined by
$a\mapsto
a\cdot \lambda.$ Clearly, $f$ is an isomorphism of $H-$modules which
satisfies
$\Delta(f(a))=(f\ot f)\tilde{\Delta}(a).$
Therefore
$(k[H],\tilde{\Delta},\varepsilon)$ is a coalgebra
where $\varepsilon=f^*(I),$ which is equivalent to
saying that $J$ satisfies (\ref{1}). In order to show that
$J$ is a twist it remains to show that it is invertible.
The invertibility of $J$ is proved in [M, Proposition 13].
We reproduce the proof (in a slighly expanded form) for the
convenience of the reader. Suppose on the contrary that $J$
is not invertible. Then there exists $0\ne L\in End_kV\otimes 
End_kV$ such that $JL=0.$ Let $F:(End_kV)^*\otimes
(End_kV)^*\rightarrow (End_kV)^*\otimes (End_kV)^*$ be
defined by $F(x)=xL.$ Clearly, $F$ is a morphism of $H\times
H$ representations, and $F\circ\tilde{\Delta}=0.$ Thus the
image $Im(F^*)$ of the morphism of $H\times H$
representations $F^*:End_kV\otimes End_kV\rightarrow
End_kV\otimes End_kV$ is contained in the kernel of the
multiplication map $m:=\tilde{\Delta}^*.$ Let
$U:=(End_kV\otimes 1)Im(F^*)(1\otimes  
End_kV).$ Clearly, $U$ is contained in the kernel of $m$ too.
But, for any $x\in U$ and $h\in H,$ $(1\otimes X_h)x(1\otimes
X_h)^{-1}\in U$ (as $End_kV=sp\{X_h|h\in H\}$ with the
relations $X_hX_{h'}=c(h,h')X_{hh'}$, i.e. conjugation by
$X_h$ is the same as the action of $h$). Thus, $U$ is a left
$End_kV\otimes End_kV$ module under left multiplication.
Similarly, it is a right module over this algebra under right
multiplication. So, it is a bimodule over $End_kV\otimes
End_kV$. Since $U\ne 0$, this implies that $U=End_kV\otimes
End_kV.$ This is a contradiction, since we get that $m=0$.
Hence $J$ is invertible.

It is straightforward to see that two different choices of $\lambda$
produce two gauge
equivalent twists $J$ for $k[H],$ so the equivalence class of $J$ is
canonically
associated to $(H,V)$.

Now, for any group $G\supseteq H,$ whose order is prime to the
characteristic of $k,$ define a triangular semisimple
Hopf algebra $F(G,H,V)=(k[G]^J,J_{21}^{-1}J)$. We wish to show that it is
also cosemisimple.
\begin{Lemma}\label{invol}
The Drinfeld element of the triangular semisimple
Hopf algebra $(A,R)=F(G,H,V)$ equals $1.$
\end{Lemma}
\proof The Drinfeld element $u$ is a grouplike element of $A,$
and for any finite-dimensional $A$-module $V$ one has
$tr|_V(u)=dim_{Rep(A)}V=dim V$ (since $Rep(A)$ is equivalent to $Rep(G),$
see [EG1, Section 1]). In particular, we can set $V$
to be the regular representation,
and find that $tr|_A(u)=dim(A)\ne 0$ in $k$. But it is clear that if
$g$ is a non-trivial grouplike element in any finite-dimensional Hopf
algebra $A,$ then $tr|_A(g)=0$. Thus, $u=1.$ \qed
\begin{Corollary} The triangular semisimple Hopf algebra $(A,R)=F(G,H,V)$
is cosemisimple.
\end{Corollary}
\proof Since $u=1,$ one has $S^2=I$ and hence $A$ is cosemisimple
(as $dim(A)\ne 0$). \qed

Thus we have assigned a triangular semisimple and cosemisimple Hopf
algebra with Drinfeld element $u=1$ to any triple $(G,H,V)$ as above.
\section{The Classification in Characteristic $0$}
In this section we assume that $k$ is of characteristic $0.$
\begin{Theorem}\label{main}
The assignment $F:(G,H,V)\mapsto (A,R)$
is a bijection between isomorphism classes of triples $(G,H,V)$
where $G$ is
a finite group, $H$ is a subgroup of $G$, and
$V$ is an irreducible projective representation of $H$
over $k$ satisfying $dimV=|H|^{1/2}$, and isomorphism classes of
triangular semisimple Hopf algebras over $k$ with Drinfeld element $u=1$.
\end{Theorem}
\proof We need to construct an assignment $F'$ in the other direction, and
check that both $F'\circ F$ and $F\circ F'$ are the identity assignments.

Let $(A,R)$ be a triangular
semisimple Hopf algebra over $k$ whose
Drinfeld element $u$ is $1.$
It follows from Deligne's theorem on Tannakian categories (see
[EG1, Theorem 2.1]) that there exist finite groups $H\subseteq G,$
and a {\em minimal} twist $J\in k[H]\otimes
k[H]$, such that $(A,R)\cong (k[G]^J,J_{21}^{-1}J)$ as triangular Hopf
algebras. As we proved in Section 3, these data are unique 
up to isomorphism and gauge transformations. 

Following Movshev
[M], define a coalgebra $B_J$ which is $k[H]$ as a vector space,
with coproduct $\tilde{\Delta}(x)=(x\otimes x)J$, $x\in H$, and the usual
counit.
This coalgebra has a natural $H$-action by left multiplication.
It follows from [M] that the coalgebra $B_J$ is simple (see [EG3] for more
explanations). Thus, the dual algebra $B_J^*$ is simple as well, and carries
an action of $H$. So we see that
$B_J^*$ is isomorphic to $End_{k}(V)$ for some vector space $V$,
and we have a homomorphism $\pi:H\to PGL(V)$.
Thus $V$ is a projective representation of $H$. It is shown in [M,
Proposition 8] that
this representation is irreducible, and it is obvious that
$\text{dim}V=|H|^{1/2}$.

It is clear that the isomorphism class
of the representation $V$ does not change if $J$ is replaced by a
twist $J'$ which is gauge
equivalent to $J$ as twists for $k[H].$
Thus, to any isomorphism class of triangular semisimple Hopf algebras
$(A,R)$ over $k$ with
Drinfeld element 1, we have assigned an isomorphism class of triples
$(G,H,V)$. Let us write this as $(G,H,V)=F'(A,R)$.

Thus, we have constructed the map $F'$.

The identity $F\circ F'=id$ follows
from [M, Proposition 5]. The identity $F'\circ F=id$
follows from Theorem \ref{uni}.
\qed

Now let $(G,H,V,u)$ be a quadruple, in which $(G,H,V)$ is as above,
and $u$ is a central element of $G$ of order $\le 2$.
We extend the map $F$ to quadruples by setting
$F(G,H,V,u)=(A,RR_u)$, where $(A,R)=F(G,H,V)$.
\begin{Theorem}\label{clas}
The assignment $F$ is a bijection between isomorphism classes of
quadruples $(G,H,V,u)$ where $G$
is a finite group,
$H$ is a subgroup of $G$, $V$ is an irreducible
projective representation of $H$ over $k$ satisfying $dimV=|H|^{1/2},$ and $u\in G$ is
a central element of order $\le 2,$
and isomorphism classes of triangular
semisimple Hopf algebras over $k.$
\end{Theorem}
\proof Define $F'$ by $F'(A,R)=(F'(A,RR_u),u),$ where 
$F'(A,RR_u)$ is defined in the proof of Theorem 5.1.
It is straightforward to see that both $F'\circ F$ and
$F\circ F'$ are the identity assignments. \qed

Theorems \ref{main} and \ref{clas} imply the following
classification
result for minimal triangular semisimple Hopf algebras over $k.$
\begin{Proposition}\label{gen} $F(G,H,V,u)$ is minimal
if and only if $G$ is generated by $H$ and $u$.
\end{Proposition}
\proof As we have already pointed out, if $(A,R)=F(G,H,V)$
then
the sub Hopf algebra $k[H]^J\subseteq A$ is minimal triangular.
Therefore, if $u=1$ then $F(G,H,V)$ is minimal if and only if $G=H$. This
obviously
remains true for $F(G,H,V,u)$ if $u\ne 1$ but $u\in H$. If $u\notin H$
then it is clear that the $R-$matrix of $F(G,H,V,u)$ generates $k[H']$,
where $H'=H\cup uH$. This proves the proposition.
\qed

\begin{Remark} {\rm As was pointed out already by Movshev, the theory
developed
in [M] and extended here is an analogue, for finite groups, of
the theory of quantization of skew-symmetric solutions of the classical
Yang-Baxter equation, developed by Drinfeld [Dr2]. In particular,
the operation $F$ is the analogue of the operation of quantization
in [Dr2].}
\end{Remark}
\section{The Classification in Positive Characteristic}
In this section we assume that $k$ is of
positive characteristic $p$, and prove an analogue of Theorem \ref{clas}
by using this theorem itself and the lifting techniques from
[EG4].

Let $F$ be the assignment from Section 4. We now have the following.
\begin{Theorem}\label{clasp}\footnote{The proof is incomplete but a complete proof is given in \cite[Remark 4.4]{e}.}
The assignment $F$ is a bijection between isomorphism classes of
quadruples $(G,H,V,u)$ where $G$ is a finite group
of order prime to $p$, $H$ is a subgroup of
$G$, $V$ is an irreducible projective representation of $H$ over $k$ satisfying
$dimV=|H|^{1/2},$ and $u\in G$ is a central element of order $\le 2,$
and isomorphism classes of triangular semisimple and cosemisimple Hopf
algebras over $k.$
\end{Theorem}
\proof As in the proof of Theorem \ref{clas} we need to construct the
assignment $F'.$ 

We recall some notation from [EG4].
Let ${\cal O}=W(k)$ be the ring of Witt vectors of $k,$ and $K$ the field
of fractions of ${\cal O}$ (it is of characteristic $0$). Let $\bar K$ be
the algebraic closure of $K.$

Let $(A,R)$ be a triangular semisimple and cosemisimple Hopf
algebra over $k.$ Lift it (see [EG4]) to a triangular semisimple Hopf algebra
$(\bar A,\bar R)$ over $K.$ By Theorem \ref{clas} we have that $(\bar 
A\ot_{K}\bar K,\bar R)=F(G,H,V,u).$
We can now reduce $V$ ``$mod\,p$'' to get $V_p$ which is an irreducible
projective representation of $H$ over the field $k.$ This can be done
since $V$ is defined by a nondegenerate $2-$cocycle $c$  (see [M])
with values in roots of unity of degree $|H|^{1/2}$
(as the only irreducible representation of the simple $H$-algebra
with basis $\{X_h|h\in H\}$, and relations $X_gX_h=c(g,h)X_{gh}$).
This cocycle can be reduced $mod\,p$ and remains nondegenerate
(since the groups of roots of unity of order $|H|^{1/2}$ in $k$ and $K$
are naturally isomorphic),
so it defines an irreducible projective representation $V_p.$ Define
$F'(A,R)=(G,H,V_p,u).$ It is shown like in characteristic $0$ that
 $F\circ F'$ and $F'\circ F$ are the
identity assignments. \qed
\begin{Corollary}\label{byatwist}
Any triangular semisimple and
cosemisimple Hopf algebra over $k$ is obtained from a group algebra
by a twist.
\end{Corollary}
\begin{Remark} {\rm Previously the statement of the corollary was only
known in characteristic $0$ [EG1, Theorem 2.1], and the best result known
to us in positive characteristic was [EG4, Theorem 3.9].}
\end{Remark}
\begin{Proposition}\label{genp}
Proposition \ref{gen} holds in positive characteristic.
\end{Proposition}
\proof As before,
if $(A,R)=F(G,H,V)$ then
the sub Hopf algebra $k[H]^J\subseteq A$ is minimal triangular.
This follows from the facts that it is true in characteristic $0,$
and that the rank of a triangular structure
does not change under lifting. Thus, Proposition
\ref{gen} holds in characteristic $p$. \qed
\begin{Remark} {\rm In view of the results of this section, the results
of our previous paper [EG3] generalize, without changes, to
cotriangular semisimple
and cosemisimple Hopf algebras in positive characteristic.}
\end{Remark}
\section{The Solvability of the Group Underlying a Minimal Triangular
Semisimple Hopf Algebra}
A classical fact about complex representations of finite groups is that
the dimension of any irreducible representation of a finite group $K$ does
not exceed $|K:Z(K)|^{1/2}$, where $Z(K)$ is the center of $K$.
Groups of central type are those groups for which this inequality
is in fact an equality. More precisely, a finite group $K$ is said to be of
{\em central type} if it has an
irreducible representation $V$ such that $(dimV)^2=|K:Z(K)|$ (see e.g. [HI]).
We shall need the following theorem 
(conjectured by Iwahori and Matsumoto in 1964)
whose proof uses the classification of
finite simple groups.
\begin{Theorem} {\bf [HI, Theorem 7.3]}\label{ct}
Any group of central type is solvable.
\end{Theorem}

As corollaries, we have the following results.
\begin{Corollary}\label{solv}
Let $A$ be a minimal triangular semisimple
Hopf algebra over $k,$
and let $G$ be the finite group such that the categories of
representations $Rep(G)$ and $Rep(A)$ are equivalent. Then $G$ is
solvable.
\end{Corollary}
\proof
We may assume that $k$ has characteristic $0$ (otherwise we can lift to
characteristic $0$),
and by Proposition \ref{gen}, that the
Drinfeld element of $A$ is $1.$ By Theorem \ref{main}, the
corresponding group $G$ has an
irreducible projective representation $V$ with $dimV=|G|^{1/2}.$ Let $K$
be a finite central extension
of $G$ with central
subgroup $Z,$ such that $V$ lifts to a linear representation of $K$.
We have $\text{(dimV)}^2=|K:Z|$.
Since $\text{(dimV)}^2\le |K:Z(K)|$ we get that $Z=Z(K)$ and hence that
$K$ is
a group of central type. But by Theorem \ref{ct}, $K$ is solvable and
hence $G\cong K/Z(K)$ is solvable as well. \qed
\begin{Remark}\label{new} {\rm Movshev conjectures in the
introduction to [M] that any finite group with a
nondegenerate $2-$cocycle is solvable. As explained in the
Proof of Corollary \ref{solv}, this result forllows from
Theorem \ref{ct}.}
\end{Remark}
\begin{Corollary}\label{rep}
Let $A$ be a triangular semisimple and cosemisimple Hopf algebra over $k$
of dimension bigger than $1$. Then $A$ has a non-trivial grouplike
element.
\end{Corollary}
\proof
We can assume that the Drinfeld element $u$ is equal to $1$ and that $A$
is not cocommutative.
Let $A_0$ be the minimal part of $A$. By Corollary \ref{solv},
$A_0=k[H]^J$ for a solvable group $H$, $|H|>1$. Therefore,
$A_0$ has non-trivial 1-dimensional representations. Since
$A_0\cong A_0^{*op}$ as Hopf algebras, we get that $A_0$,
and hence $A$, has non-trivial grouplike elements.\qed

Corollary \ref{rep} motivates the following question.
\begin{Question}\label{ques}
Let $(A,R)$ be a quasitriangular semisimple and cosemisimple
Hopf algebra
over $k$ (e.g. the quantum double of a semisimple and
cosemisimple Hopf algebra), and let
$dim(A)>1$. Is it true that $A$ possesses a non-trivial grouplike
element?
\end{Question}
Note that a positive answer to this question would imply that
for a semisimple and cosemisimple Hopf algebra $A$ over $k$ either $A$ or
$A^*$
possesses a non-trivial grouplike element. Such a result is very
desirable for the problem of the classification of semisimple Hopf
algebras.
\section{Group-Theoretical Data Corresponding to Minimal Triangular
Hopf Algebras Constructed from a Bijective 1-Cocycle}
In this section we determine the group-theoretical data
corresponding, under the bijection of the classification given in Theorem
\ref{main},
to the minimal triangular semisimple Hopf algebras constructed in
[EG2, Section 4]. We will use the definitions and notation from [EG2].

Let $k=\C$.
Let $G$ be a finite group,
$A$ be a finite abelian group with a left $G$-action $(g,a)\mapsto g\cdot
a,$
and $\pi:G\to A$ be a bijective 1-cocycle, i.e. a bijective map
such that $\pi(gg')=\pi(g)+g\pi(g')$ (in particular, $|G|=|A|$).
Let $H=G\propto A^*$ be the semidirect product of $G$ by the dual group
$A^*$ to $A.$ Following [EG2, Section 4], we can associate
to this data the element
$$
J=|A|^{-1}\sum_{g\in G,b\in A^*}e^{(\pi(g),b) }b\otimes g
$$
(for convenience we use the opposite element to the one from [EG2]).
We proved in [EG2]
that this element is a minimal twist for $k[H],$ so $k[H]^J$ is a minimal
triangular semisimple Hopf algebra with Drinfeld element $u=1.$
Now we wish to find the irreducible projective representation $V$ of
$H$ which corresponds to $k[H]^J$ under the correspondence of 
 Theorem \ref{main}.

We will now construct the representation $V$ and show that it is the one
corresponding to $k[H]^J.$

Let $V=Fun(A,k)$ be the space of $k$-valued functions on $A$.
It has a basis $\{\delta_a|a\in A\}$ of characteristic functions of
points. Define a projective action $\phi$ of $H$ on $V$
by
$$\phi(b)\delta_a=e^{-(a,b)}\delta_a,\; \phi(g)\delta_a=\delta_{g\cdot
a+\pi(g)}\;{\rm and}\;\phi(bg)=\phi(b)\phi(g)$$
for $g\in G$ and $b\in A^*.$ It is straightforward to
verify that this is indeed a projective representation.
\begin{Proposition}\label{projrep}
The representation $V$ is irreducible, and corresponds
to $k[H]^J$ under the bijection of the classification given in Theorem
\ref{main}.
\end{Proposition}
\proof Let $B_J$ be the coalgebra $(k[H],\tilde\Delta)$ where
$\tilde\Delta(x)=(x\otimes x)J$, $x\in H$. Let $B_J^*$ be the dual
algebra.
It is enough to show that the $H$-algebras $B_J^*$ and $End_k(V)$
are isomorphic.

Let us compute the multiplication in the algebra $B_J^*$.
We have
\begin{equation}
\tilde\Delta(bg)=|A|^{-1}\sum_{g'\in G,b'\in B}
e^{(\pi(g'),b')}b(g\cdot b')g\otimes bgg'.\label{2}
\end{equation}
Let $\{Y_{bg}\}$ be the dual basis of $B_J^*$ to the basis $\{bg\}$ of
$B_J.$
Let $*$ denote the multiplication law dual to the coproduct
$\tilde\Delta$.
Then, dualizing equation (\ref{2}), we have
\begin{equation}
Y_{b_2g_2}*Y_{b_1g_1}=e^{(\pi(g_1)-\pi(g_2),b_2-b_1)}Y_{b_1g_2}\label{3}
\end{equation}
for $g_1,$ $g_2\in G$ and $b_1,$ $b_2\in A^*$
(here for convenience we write the operations in $A$ and $A^*$ additively).
Define $Z_{bg}:=e^{(\pi(g),b)}Y_{bg}$ for $g\in G$ and $b\in A^*.$ In the
basis
$\{Z_{bg}\}$ the multiplication law in $B_J^*$ is given by
\begin{equation}
Z_{b_2g_2}*Z_{b_1g_1}=e^{(\pi(g_1),b_2)}Z_{b_1g_2}.\label{4}
\end{equation}

Now let us introduce a left action of $B_J^*$ on $V$. Set
\begin{equation}
Z_{bg}\delta_a=e^{(a,b)}\delta_{\pi(g)}.\label{5}
\end{equation}
It is straightforward to check using (4) that (5) is indeed a left action.
It is also straightforward to compute that this action is $H-$equivariant.
Thus, (5) defines an isomorphism $B_J^*\to End_k(V)$ as $H-$algebras,
which proves the proposition. \qed

\end{document}